\documentclass[11pt]{amsart}
\usepackage{amsmath,amsfonts,amssymb,amsxtra,latexsym,amscd,enumerate,amsthm}

\usepackage{latexsym}
\usepackage{epsfig}

\def\eps{\epsilon}

\def\R{\mathbb{R}}
\def\C{\mathbb{C}}

\def\N{\mathcal{N}}
\def\S{\mathcal{S}}
\def\sgn{\operatorname{sgn}}
\def\Re{\operatorname{Re}}

\def\beq{\begin{equation}}
\def\eeq{\end{equation}}

\def\neq{\nLeftrightarrow}
\def\beq{\begin{equation}}
\def\eeq{\end{equation}}

\newtheorem{definition}{Definition}
\newtheorem{proposition}{Proposition}
\newtheorem{theorem}{Theorem}
\newtheorem{lemma}{Lemma}
\newtheorem{corollary}{Corollary}

\begin{document}

\title[Sharp well-posedness for NLS]{Sharp well-posedness and ill-posedness results for a quadratic non-linear Schr\"odinger equation}

\author{Ioan Bejenaru}
\address{Department of Mathematics, UCLA, Los Angeles CA 90095-1555}
\email{bejenaru@math.ucla.edu}

\author{Terence Tao}
\address{Department of Mathematics, UCLA, Los Angeles CA 90095-1555}
\email{tao@math.ucla.edu}

\begin{abstract}  We establish that the quadratic non-linear Schr\"odinger equation
$$ iu_t + u_{xx} = u^2$$
where $u: \R \times \R \to \C$, is locally well-posed in $H^s(\R)$ when $s \geq -1$ and ill-posed when $s < -1$.  Previous work in \cite{kpv} had established local well-posedness for $s > -3/4$.  The local well-posedness is achieved by an iteration using a modification of the standard $X^{s,b}$ spaces.  The ill-posedness uses an abstract and general argument relying on the high-to-low frequency cascade present in the non-linearity, and a computation of the first non-linear iterate.
\end{abstract}

\maketitle

\section{Introduction}

Consider the Cauchy problem 
\begin{equation}\label{idelta}
\begin{split}
 i u_t + \Delta u &= F(u) \\
 u(0) &= f \in H^s_x(\R^n) \\
 u &\in C^0_t H^s_x([0,T] \times \R^n)
\end{split}
\end{equation}
for a semilinear Schr\"odinger equation on $[0,T] \times \R^n$ for some local time interval\footnote{One can also consider the Cauchy problem
backwards in time, on some interval $[-T,0]$, but this backwards problem is equivalent to the forwards problem after applying the conjugation
$u(t) \mapsto \overline{u(-t)}$, $F(z) \mapsto \overline{F(\overline{z})}$.} $[0,T]$ with $T > 0$, where
the initial data $f$ is given in some Sobolev space\footnote{We will subscript spatial function spaces by $x$ and temporal function spaces by $t$, thus for instance $C^0_t H^s_x([0,T] \times \R^n)$ is the space of all functions $u(t,x)$ for which the map $t \mapsto u(t)$ is continuous into $H^s_x(\R^n)$, equipped with the norm $\sup_{0 \leq t \leq T} \|u(t) \|_{H^s_x(\R^n)}$.}  $H^s_x(\R^n)$, the solution $u$ is complex-valued,
 and $F: \C \to \C$ is a power-type nonlinearity (thus $|F(z)| \sim |z|^p$ for some exponent $p$, and similarly for derivatives).
 To fix conventions, we define the Sobolev space $H^s_x(\R^n)$ for any $s \in \R$ as the Banach space of distributions $f$ for which the norm
$$ \| f \|_{H^s_x(\R^n)} := \| \langle \xi \rangle^s \hat f \|_{L^2_\xi(\R^n)}$$
is finite, where $\langle \xi \rangle := (1 + |\xi|^2)^{1/2}$, and $\hat f$ is the Fourier transform
$$ \hat f(\xi) := \int_{\R^n} e^{ix \cdot \xi} f(x)\ dx.$$
 This particular problem has been studied extensively in the literature, for various values of $n$, $s$, and $F$, as it is a simple
 model for the more general Cauchy problem for non-linear dispersive equations.  
In the situation considered in this paper, the regularity $H^s_x(\R^n)$ is very low (in fact, $s$ will be negative), so that the solutions to \eqref{idelta} cannot be interpreted in the classical sense; we will make sense of the equation for rough data later, but suffice to
say for now that we will be able to show that the rough $C^0_t H^s_x$ solutions we construct will be strong limits in $C^0_t H^s_x$ 
of smooth solutions.

If $F$ is smooth, then one typically obtains a local well-posedness result\footnote{By this we mean that for any choice of initial data $u_0 \in H^s_x(\R^n)$ there exists a time $T > 0$, and a continuous solution map defined in a small ball in $H^s_x(\R^n)$ centered at $u_0$, and taking values in $C^0_t H^s_x([-T,T] \times \R^n)$.  Furthermore, when the data is restricted to a suitable smooth class (e.g. $H^s(\R^n)$ for $s>n/2$), then the solution map agrees with the standard (and unique) solutions that can be constructed for instance by the energy method.}
 when $s$ is large, but not when $s$ is small.  For instance, for the power-type semilinear equation
\begin{equation}\label{nls-power}
 i u_t + \Delta u = \pm |u|^{p-1} u
\end{equation}
for $p > 1$ and either choice of sign $\pm$, this equation is locally well-posed\footnote{If $p$ is not an odd integer, then we also need
the technical condition $p > \lfloor s \rfloor + 1$ to ensure that the non-linearity is at least as regular as the nonlinearity.}
 when $s \geq \min(0, s_c)$, where the scaling regularity $s_c$ is defined by
$$ s_c = \frac{n}{2} - \frac{2}{p-1};$$
see for instance \cite{cwi}.  This condition is fairly sharp; if $s < \min(0,s_c)$, then the solution map is known to not be uniformly continuous
from $H^s$ to $C_t H^s([-T,T] \times \R^n)$, even if we make $T$ small and restrict the data to a small ball around the origin; see \cite{bkpv}, \cite{kpv-ill}, \cite{cct}.
The regularity $H^{s_c}(\R^n)$ is a natural limit to well-posedness as it is preserved by the scale invariance
\begin{equation}\label{scaling}
 u(t,x) \mapsto \lambda^{-2/(p-1)} u(\frac{t}{\lambda^2}, \frac{x}{\lambda}) \hbox{ for any } \lambda > 0
 \end{equation}
of the equation \eqref{nls-power}, while the 
regularity $L^2(\R^n)$ is another natural limit, as it is preserved by the Galilean invariance
$$ u(t,x) \mapsto \exp( i (v \cdot x - |v|^2 t ) ) u( t, x - vt ) \hbox{ for any } v \in \R^n$$
of the same equation.

Thus it would seem that the local well-posedness theory for semilinear Schr\"odinger equations with power-type nonlinearity is complete.  However, it was observed in \cite{kpv} that one can lower the regularity threshold for local well-posedness below $s=0$ (i.e. below the Galilean threshold) by choosing a non-linearity which is not Galilean invariant.  In particular, the one-dimensional quadratic semilinear
Schr\"odinger equation\footnote{We shall only consider scalar solutions here for simplicity.  However, one can extend the analysis here to finite-dimensional systems with a quadratic form nonlinearity $Q(u,u)$ which is linear (as opposed to anti-linear) in both variables, without any difficulty.  Also there is no distinction between the $+u^2$ and $-u^2$ non-linearities, as can be seen by the transformation $u \mapsto -u$.}
\begin{equation}\label{nls-quad}
\begin{split}
i u_t + u_{xx} &= u^2\\
 u(0) &= f \in H^s_x(\R) \\
 u &\in C^0_t H^s_x([0,T] \times \R)
\end{split}
\end{equation}
was shown in \cite{kpv} to be locally well-posed in $H^s_x(\R)$ for all $s > -3/4$, by means of an iteration argument in 
the $X^{s,b}$ spaces; in contrast, with a quadratic non-linearity such as $|u| u$, the lowest Sobolev regularity for which one has well-posedness is $L^2_x(\R)$ (see \cite{tsutsumi}, \cite{cwi}, \cite{bkpv}, \cite{kpv-ill}, \cite{cct}).
One should remark that these regularities are well above the scaling regularity, which in this case is $s_c = -3/2$; thus these results are subcritical with respect to scaling.

The paper \cite{kpv} also considered other quadratic non-linearities such as $u \overline{u}$ and $\overline{u}^2$, obtaining similar results.  However we wish to focus on the $u^2$ nonlinearity in \eqref{nls-quad} to point out one interesting feature of this equation, namely its complex analyticity in $u$.  This manifests itself in a number of ways; in particular, this equation in the spacetime-frequency domain $(\tau,\xi)$ is almost entirely supported in the upper half-space $\tau > 0$.  To see this heuristically, let us formally introduce the spacetime Fourier transform
$$ \tilde u(\tau,\xi) := \int_\R \int_\R u(t,x) e^{i(t \tau + x \xi)}\ dt dx$$
(ignoring for now the issue of extending $u$ globally in time), then \eqref{nls-quad} transforms (heuristically, at least)
to the integral equation
\begin{equation}\label{tu}
\tilde u(\tau,\xi) = \delta(\tau - \xi^2) \hat u_0(\xi) + \frac{1}{\tau - \xi^2}
\int\int_{\tau = \tau_1 + \tau_2; \xi = \xi_1 + \xi_2} \tilde u(\tau_1, \xi_1) \tilde u(\tau_2, \xi_2)
\end{equation}
where $\delta$ is the Dirac delta function, and we will be deliberately vague about how to define in a distributional sense the operation of dividing by $\tau - \xi^2$.  
If one tries to solve the equation \eqref{tu} iteratively, viewing this equation as a way to obtain a new approximation of $\tilde u$ from an old one, starting from (say) the zero solution $\tilde u = 0$, we see that all the iterates $\tilde u$ are supported on the upper half-plane $\tau > 0$, and thus we expect the final solution to do so also.

This additional property of the problem \eqref{nls-quad} suggests that perhaps some further improvement to $H^{-3/4}_x(\R)$ and beyond is possible;
for instance, \eqref{tu} suggests that the solution is unlikely concentrate near the spacetime frequency 
origin $(\tau,\xi) = (0,0)$.  It is also similarly difficult for the iterates of the solution to return back to the parabola $\tau = \xi^2$,
where the solution is expected to concentrate (in analogy with the linear solution).  

A first step in this direction was made by Muramutu and Taoka \cite{mur}, obtaining local well-posedness in the Besov space $B^{-3/4,1}_2(\R)$, which 
is slightly stronger than $H^{-3/4}_x(\R)$, by a refinement of the $X^{s,b}$ iteration method.  However, it was shown in \cite{kpv}, \cite{nak} that the key bilinear
$X^{s,b}$ estimate needed to apply this method failed for $H^s_x(\R)$ for any $s \leq -3/4$.  Nevertheless, it turns out that we can use the additional
information that $\tilde u$ concentrates on the upper half-plane to avoid most of the counterexamples in \cite{kpv}, \cite{nak}, and after modification of the $X^{s,b}$ spaces we can in fact avoid the remaining counterexamples also, to obtain our first main theorem:

\begin{theorem}[Local well-posedness in $H^{-1}_x(\R)$]\label{lwp}  Let $r > 0$ be any radius, and let $B_r$ be the ball
$$ B_r := B_{H^{-1}_x(\R)}(0,r) := \{ u_0 \in H^{-1}_x(\R): \| u_0 \|_{H^{-1}_x(\R)} < r \}.$$
Then there exists a time $T > 0$ (in fact we obtain $T = \max(1, c r^{-1/2} )$ for some absolute constant $c>0$)
and a map $f \mapsto u[f]$ which is continuous from $B_r$ to $C^0_t H^{-1}_x([0,T] \times \R)$, such that the
restriction of this map to $B_r \cap H^s_x(\R)$ (with the $H^s_x(\R)$ topology) maps continuously 
to $C^0_t H^s([0,T] \times \R)$ for any $s \geq -1$.  Furthermore, if $f$ lies in
a smooth space, say $B_r \cap H^3_x(\R)$, then $u[f]$ lies in $C^0_t H^3_x \cap C^1_t H^1_x([0,T] \times \R)$ and
solves the equation \eqref{nls-quad} in the classical sense.   
\end{theorem}

We prove this theorem in Section \ref{lwp-sec}, as a consequence of standard iteration machinery and a construction of a function space
obeying certain linear and bilinear estimates.
It is easy to establish (e.g. by energy methods) that classical solutions to \eqref{nls-quad} in $C^0_t H^3_x \cap C^1_t H^1_x$ are unique, and so 
the solution map $S$ given by the above theorem is the unique strong limit of smooth solutions in $C^0_t H^{-1}_x$.  The theorem also shows
that if the solution blows up at some time $T_*$, then the $H^{-1}_x(\R)$ norm of $u(t)$ must blow up at a rate of $c|T_*-t|^{-2}$ or greater
as $t \to T_*$.  The main novelty in the proof is a modification of the $X^{s,b}$ spaces in order to exploit the concentration of the solution in the upper half-plane $\tau > 0$, and also to deal with the failure
of the $X^{s,b}$ norms to adequately control the behavior of this equation near the time axis $\xi = 0$.

Our second main result is that the threshold $s \geq -1$ is completely sharp.

\begin{theorem}[Ill-posedness below $H^{-1}_x(\R)$]\label{ill-p}  Let $r > 0$ be arbitrary, and let $T$ and $f \mapsto u[f]$ be as in Theorem \ref{lwp}.
Then the solution map $f \mapsto u[f]$ is discontinuous on $B_r$ (with the $H^s_x(\R)$ topology) to $C^0_t H^{-1}_x([0,T] \times \R)$ 
(with the $C^0_t H^{s'}_x([0,T] \times \R)$ topology) for
any $s <-1$ and $s'\in \R$.
\end{theorem}

This theorem will be proven by demonstrating a high-to-low frequency cascade in the first non-trivial iterate of an integral equation 
associated to \eqref{nls-quad}; see Section \ref{lwp-sec}.  We will then invoke a rather general result (Proposition \ref{power} below), 
which may be of independent interest,
which shows that any non-linear evolution equation with polynomial nonlinearity will be illposed whenever a high-to-low frequency cascade
in one of its iterates can be established.

The authors thank Martin Hadac and Sebastian Herr for pointing out an error in the original manuscript, and Hideo Takaoka and Zhenqiu Zhang for further corrections.  The second author is supported by a grant from
the Packard Foundation.

\section{Reduction to an integral equation}\label{reduced-sec}

We first give some very standard reductions for Theorem \ref{lwp}.  The first is to use the scale invariance \eqref{scaling} to
scale the radius $r$ to be small.  Indeed
if one defines
$$ f^{(\lambda)}(x) := \frac{1}{\lambda^2} f(\frac{x}{\lambda})$$
then a simple computation shows\footnote{Here and in the sequel we use $C, c > 0$ to denote various positive absolute constants.}
$$ \| f^{(\lambda)} \|_{H^{-1}_x(\R)} \leq C \lambda^{-1/2} \| f \|_{H^{-1}_x(\R)}$$
for $\lambda > 1$.  We can thus rescale the initial data to be arbitrarily small in
$H^s_x(\R)$ norm.  It thus suffices to prove Theorem \ref{lwp} when $T = 1$ and $r$ is sufficiently small (smaller than some absolute
constant $c > 0$).

Next, we shall use Duhamel's formula to recast \eqref{nls-quad} in the integral form
\begin{equation}\label{duhamel}
u(t) = \exp(it\partial_{xx}) f + \int_0^t \exp(i(t-s)\partial_{xx})(u(s)^2)\ ds,
\end{equation}
where $\exp(it\partial_{xx})$ is the propagator for the free Schr\"odinger equation $iu_t + u_{xx} = 0$, which can be defined for instance using the
spatial Fourier transform as
$$ \widehat{\exp(it\partial_{xx}) f}(\xi) := e^{-it \xi^2} \hat f(\xi).$$
Following Bourgain \cite{borg:xsb}, it turns out to be convenient to replace the local-in-time integral equation \eqref{duhamel} with a global-in-time
truncated integral equation.  Let $\eta: \R \to \R$ be a smooth bump function such that $\eta(t) = 1$ for $|t| \leq 1$ and $\eta(t) = 0$ for $|t| > 2$, and let $a(t) := \frac{1}{2} \sgn(t) \eta(t/5)$.  
Then observe that
$$ \int_0^t g(s)\ ds = \eta(t) \int_\R a(s) g(s)\ ds + \int_\R a(t-s) g(s)\ ds$$ 
for all $0 \leq t \leq 1$ and any $g: \R \to \R$.  Hence we can replace \eqref{duhamel} on the time interval $0 \leq t \leq 1$ by the equation
\begin{equation}\label{duhamel-trunc}
\begin{split}
u(t) &= \eta(t) \exp(it\partial_{xx}) f + \eta(t) \exp(it\partial_{xx}) \int_\R a(s) \exp(-is\partial_{xx})(u(s)^2)\ ds\\
&\quad + \int_\R a(t-s) \exp(i(t-s)\partial_{xx})(u(s)^2)\ ds,
\end{split}
\end{equation}
in the sense that any classical (e.g. $C^0_t H^3_x(\R \times \R)$ will do) global-in-time solution to \eqref{duhamel-trunc} is also a classical solution to \eqref{duhamel} and hence \eqref{nls-quad}.  Note that if $u \in C^0_t H^3_x$, one can easily use \eqref{duhamel}
or \eqref{duhamel-trunc} and Sobolev embedding to conclude that $u \in C^1_t H^1_x$, and so one can make sense of \eqref{nls-quad} in a classical sense.

It remains to find global-in-time solutions to \eqref{duhamel-trunc} for initial data $f$ in $B_r$.

We will write \eqref{duhamel-trunc} more abstractly as
\begin{equation}\label{uln-2}
 u = L(f) + N_2(u,u)
 \end{equation}
where $L$ is the linear operator
\begin{equation}\label{Ldef}
L(f)(t) := \eta(t) \exp(it\partial_{xx}) f
\end{equation}
and $N_2$ is the bilinear operator
\begin{equation}\label{Ndef}
\begin{split}
N_2(u,v)(t) &:=
\eta(t) \exp(it\partial_{xx}) \int_\R a(s) \exp(-is\partial_{xx})(u(s)v(s))\ ds\\
&+ \int_\R a(t-s) \exp(i(t-s)\partial_{xx})(u(s) v(s))\ ds.
\end{split}
\end{equation}
The subscript $2$ denotes the fact that this operator is quadratic.
We now pause to systematically develop the well-posedness, persistence of regularity, 
and ill-posedness theory for such an abstract type of operator.
This theory is mostly standard, but the material on ill-posedness may be of independent interest.

\section{Abstract well-posedness and ill-posedness theory}

In this section we shall consider an abstract semilinear evolution equation with a $k$-linear non-linearity for some $k \geq 2$.  Specifically,
we consider the abstract equation
\begin{equation}\label{uln}
u = L(f) + N_k(u, \ldots, u)
\end{equation}
where the initial data $f$ takes values in some data space $D$,
the solution $u$ takes values in some solution space $S$, 
the linear operator $L: D \to S$ is densely defined, and the $k$-linear
operator $N_k: S \times \ldots \times S \to S$ is also densely defined.  

\begin{definition}[Quantitative well-posedness]  Let $(D, \| \|_D)$ be a Banach space of initial data, and $(S, \| \|_S)$ be
a Banach space of spacetime functions. We say that the equation \eqref{uln} is \emph{quantitatively well posed} in
the spaces $D,S$ if one has estimates of the form\footnote{We adopt the convention that if $X$ is a Banach space, then $\|u\|_X$ denotes the norm of $u$ in $X$, and that $\|u\|_X = \infty$ if $u \not \in X$.}
$$ \| L(f) \|_S \leq C \| f \|_D$$
and
$$ \| N_k( u_1, \ldots, u_k ) \|_S \leq C \| u_1 \|_S \ldots \| u_k \|_S$$
for all $f \in D$, $u_1,\ldots,u_k \in S$ and some constant $C > 0$.  
\end{definition}

Note that once the equation \eqref{uln} is quantitatively well-posed, one can extend the densely defined operators $L$ and $N_k$ to all of $D$ and $S \times \ldots \times S$ respectively in the unique continuous fashion.

If $(D, \| \|_D)$ is a Banach space, we use $B_D(0,r) := \{ f \in D: \|f\|_D < r \}$ to denote the usual open ball of radius $r$ around
the origin.  The standard well-posedness result for such equations is that quantitative well-posedness implies analytic well-posedness.
More precisely:

\begin{theorem}[Standard well-posedness theorem]\label{swp}  Suppose the equation \eqref{uln} is quantitatively well posed
in the spaces $D,S$.  Then there exists constants $C_0,\eps_0 > 0$ such that for 
all $f \in B_D(0,\eps_0)$, there there exists a unique solution $u[f] \in B_S(0, C_0 \eps_0)$ to the equation 
\eqref{uln}.  More specifically, if we define the nonlinear maps $A_n: D \to S$ for $n=1,2,\ldots$ by the recursive formulae
\begin{align*}
A_1(f) &:= L(f) \\
A_n(f) &:= \sum_{n_1,\ldots,n_k \geq 1: n_1 + \ldots + n_k = n} N_k( A_{n_1}(f), \ldots, A_{n_k}(f) ) \hbox{ for } n > 1
\end{align*}
then we have the homogeneity property
\begin{equation}\label{homog-est}
A_n(\lambda f) = \lambda^n A_n(f) \hbox{ for all } \lambda \in \R, n \geq 1, \hbox{ and } f \in D
\end{equation}
(so in particular $A_n(0) = 0$) and the Lipschitz bound
\begin{equation}\label{lip}
 \| A_n(f) - A_n(g) \|_S \leq \|f-g\|_D C^n_1 (\|f\|_D + \|g\|_D)^{n-1}
\end{equation}
for some $C_1 > 0$, all $f,g \in D$, and all $n \geq 1$.  In particular we have
\begin{equation}\label{lip-0}
 \| A_n(f) \|_S \leq C^n_2 \|f\|_D^n
\end{equation}
for some $C_2 > 0$.  Furthemore, we have the absolutely convergent (in $S$) power series expansion
\begin{equation}\label{uanf}
u[f] = \sum_{n=1}^\infty A_n(f)
\end{equation}
for all $f \in B_D(0,\eps_0)$.
\end{theorem}

Thus for instance, if $k=2$, then
\begin{align*}
A_1(f) &= L(f) \\
A_2(f) &= N_2(Lf, Lf) \\
A_3(f) &= N_2(Lf, N_2(Lf,Lf)) + N_2(N_2(Lf,Lf),Lf) \\
&\vdots
\end{align*}
whereas if $k=3$, then
\begin{align*}
A_1(f) &= L(f) \\
A_2(f) &= 0 \\
A_3(f) &= N_3(Lf, Lf, Lf) \\
A_4(f) &= 0 \\
&\vdots
\end{align*}
In general, one can express $A_n(f)$ as a sum over $k$-ary trees with $n$ nodes, but we will not need such an explicit representation here.

\begin{proof}  We shall be somewhat brief here since this theorem is well-known.
For any fixed $f \in B_D(0,\eps_0)$, one can easily verify from the quantitative well-posedness
hypothesis that the map
$u \mapsto L(f) + N_k(u,\ldots,u)$ will be a contraction from $B_S(0,C_0\eps_0)$ to $B_S(0,C_0\eps_0)$ if $C_0$ is sufficiently large
and $\eps_0$ is sufficiently small (depending on $C_0$).  The contraction mapping theorem then gives the existence and uniqueness of  
the map $f \mapsto u[f]$.

Now we start proving the power series expansion.  One can easily verify \eqref{homog-est} by induction; in fact we easily check that
$$ A_n(f) = M_n(f,\ldots,f)$$
for some $n$-linear map $M_n: D \times \ldots \times D \to S$.
One can also inductively
obtain an estimate of the form
\begin{equation}\label{anfs}
\| A_n(f) \|_S \leq (C_3 \|f\|_D)^n 
\end{equation}
for some large constant $C_3 > 0$ (depending of course on the quantitative well-posedness constants and on $k$); note this
already gives \eqref{lip-0}.  We remark that for the purposes of proving  \eqref{anfs}, it is actually slightly
easier for inductive purposes to establish the slight stonger upper bound of $(C_3 \|f\|_D)^n / (C_4 n^{C_5})$, where $C_4$ and $C_5$ are somewhat 
large but not as large as $C_3$.

Now we prove \eqref{lip}.
By symmetry we may take $\|f\|_D \leq \|g\|_D$, and by scaling we can take $\|g\|_D \leq 1$.  We can of course assume that $f \not = g$.
Write $t := \|f-g\|_D$, and write $f = g+th$, thus $0 < t \leq 2$ and $\|h\|_D = 1$.  It then suffices to show 
\begin{equation}\label{anlip}
 \| A_n(g+th) - A_n(g) \|_S \leq t C^n_1.
 \end{equation}
The nonlinear operator $A_n(f)$ can be written as $M_n(f,\ldots,f)$ for some $n$-linear operator $M_n$, which implies
that for fixed $g,h$, the function $s \mapsto A_n(g+sh) - A_n(g)$ is a polynomial of degree at most $n$ in $s$ with zero constant term, thus
\begin{equation}\label{angsh}
A_n(g+sh) - A_n(g) = \sum_{j=1}^n F_j s^j
\end{equation}
for some $F_j \in S$.  From \eqref{lip-0} and the triangle inequality we also have
$$ \| A_n(g+sh) - A_n(g) \|_S \leq (4C_2)^n$$
for all $|s| \leq 1$ (say) and some constant $C > 0$.  Thus we have
$$ \| \sum_{j=1}^n F_j s^j \|_S \leq (4C_2)^n$$
for all $|s| \leq 1$.  Using the Lagrange interpolation formula to recover $F_1,\ldots,F_n$ from various sample points of
$\sum_{j=1}^n F_j s^j$, we conclude that
$$ \| F_j \|_S \leq C^n$$
for all $1 \leq j \leq n$ and some $C > 0$.  Inserting this back into \eqref{angsh} we obtain \eqref{anlip}.  (Note that if we allow $s$ to be complex, one could
also proceed using the Cauchy integral formula instead of the Lagrange interpolation formula.)

From \eqref{lip-0} we see that the series $\sum_{n=1}^\infty A_n(f)$ is absolutely convergent in $S$ for $\eps$ small enough.  If for any integer
$K \geq 1$, we let $u_K$
be the partial sum $u_K := \sum_{n=1}^K A_n(f)$, one can easily verify a formula of the form
$$ L(f) + N_k(u_K,\ldots,u_K) = u_K + \sum_{K < n \leq nK} A_{n,K}(f)$$
where $A_{n,K}(f)$ is a non-linear expression of $f$ which consists of some subset of the terms used to form $A_n(f)$.  One can then again use induction to obtain estimates of the form
$$ \| A_{n,K}(f) \|_S \leq (C \|f\|_S)^n$$
for some $C > 0$, and hence we see that if $f \in B_D(0,\eps_0)$ for $\eps_0$ sufficiently small, that
$$ \| L(f) + N_k(u_K,\ldots,u_K) - u_K \|_S \leq (C \eps_0)^n.$$
Using the contraction mapping principle again, we see that $u_K$ converges to $u[f]$ in $S$ norm, and we obtain \eqref{uanf}.
\end{proof}

From \eqref{uanf}, \eqref{lip} one can verify that the map $f \mapsto u[f]$ is continuous (in fact Lipschitz continuous) from
$B_D(0,\eps_0)$ to $B_S(0,C_0\eps_0)$ for $\eps_0$ small enough.  In fact this Lipschitz continuity can also be read off directly from
the contraction mapping theorem.

Now, we investigate continuity in both finer and coarser topologies.  The basic result for finer topologies is as follows.

\begin{theorem}[Standard persistence of regularity theorem]\label{persist}  Suppose the equation \eqref{uln} is quantitatively well posed
in the spaces $D,S$, and let $f \mapsto u[f]$ be the solution map from $B_D(0,\eps_0)$ to $u[f] \in B_S(0, C_0 \eps_0)$
constructed in Theorem \ref{swp}.  Suppose we are given spaces $(D', \|\|_{D'})$ and $(S',\|\|_{S'})$ obeying
the estimates
$$ \| L(f) \|_{S'} \leq C \| f \|_{D'}$$
and
$$ \| N_k( u_1, \ldots, u_k ) \|_{S'} \leq C \sum_{j=1}^k \| u_j \|_{S'} \prod_{1 \leq i \leq k; i \neq j} \| u_i \|_S.$$
Then, if $\eps_0$ is sufficiently small, the solution map is also continuous from $B_D(0,\eps_0) \cap D'$ (in the $D'$ topology)
to $B_S(0, C_0 \eps_0) \cap S'$ (in the $S'$ topology).
\end{theorem}

\begin{proof}  If $f \in B_D(0,\eps_0) \cap D'$ for a suitably small $\eps_0$, then the above estimates easily imply that
$u \mapsto L(f) + N_k(u,\ldots,u)$ will be a contraction in the $S'$ norm 
from $B_S(0,C_0\eps_0) \cap S'$ to $B_S(0,C_0\eps_0) \cap S'$, and the claim follows. 
\end{proof}

Now we turn to coarser topologies. The basic result here is that if the map $f \mapsto \sum_{n=1}^\infty A_n(f)$ is continuous
in a coarse topology, then each component $f \mapsto A_n(f)$ of the series is also continuous in this coarse topology.

\begin{proposition}\label{power}  Suppose that the equation \eqref{uln} is quantitatively well-posed in the Banach spaces $D$ and $S$, with a solution map
$f \mapsto u[f]$ from a ball $B_D$ in $D$ to a ball $B_S$ in $S$.
Suppose that these spaces are then given other norms $D'$ and $S'$, which are weaker than $D$ and $S$ in the sense that
$$ \|f\|_{D'} \leq C \|f\|_D, \quad \| u \|_{S'} \leq C \| u \|_S$$
for some absolute constant $C$.  (Note that $D$ is unlikely to be complete in the $D'$ norm, and similarly for $S$ and $S'$.)  Suppose that the
solution map $f \mapsto u[f]$ is continuous from $(B_D,\| \|_{D'})$ (i.e. the ball $B_D$ equipped with the $D'$ topology) to
$(B_S,\| \|_{S'})$.  Then for each $n$, the non-linear operator $A_n: D \to S$ is continuous from $(B_D, \| \|_{D'})$ to $(S, \|\|_{S'})$.
\end{proposition}

\begin{proof}  We induct on $n$, assuming that for all $n' < n$ the operator $A_{n'}: D \to S$ has already been shown to be continuous
from $(D, \| \|_{D'})$ to $(S, \| \|_{S'})$.

Let $f_m$ be a sequence in $B_D$ which converges to $f \in B_D$ in the $D'$ topology, thus $\|f_m - f \|_{D'} \to 0$.  Our task is
to show that $\| A_n(f_m) - A_n(f) \|_{S'} \to 0$.

Now let $0 < \lambda \leq 1$ be a small number to be chosen later.  By hypothesis, the map $f \mapsto u[f]$ is continuous from $(B_D, \| \|_{D'})$ to
$(B_S, \| \|_{S'})$, and hence 
$$ \lim_{m \to \infty} \| u[\lambda f_m] - u[\lambda f] \|_{S'} = 0.$$
Expanding out the power series and using homogeneity, we have
$$ \lim_{m \to \infty} \| \sum_{n'=1}^\infty \lambda^{n'} (A_{n'}(f_m) - A_{n'}(f)) \|_{S'} = 0.$$
By the induction hypothesis we already have
$$ \lim_{m \to \infty} \| \sum_{n' < n} \lambda^{n'} (A_{n'}(f_m) - A_{n'}(f)) \|_{S'} = 0$$
so on subtracting and then dividing by $\lambda^n$, we conclude that
$$ \lim_{m \to \infty} \| \sum_{n' \geq n} \lambda^{n'-n} (A_{n'}(f_m) - A_{n'}(f)) \|_{S'} = 0$$
and hence by the triangle inequality
$$ \limsup_{m \to \infty} \| A_{n}(f_m) - A_{n}(f)) \|_{S'} 
\leq \sum_{n' > n} \lambda^{n'-n} \sup_m \| A_{n'}(f_m) - A_{n'}(f) \|_{S'}.
$$
Using \eqref{lip}, we conclude
$$ \limsup_{m \to \infty} \| A_{n}(f_m) - A_{n}(f)) \|_{S'} 
\leq \sum_{n' > n} \lambda^{n'-n} (C \eps)^{n'}.
$$
The right-hand side is convergent for $\lambda$ small enough.  Taking $\lambda \to 0$, we conclude
$$ \limsup_{m \to \infty} \| A_{n}(f_m) - A_{n}(f)) \|_{S'} = 0,$$
and the claim follows.
\end{proof}

This proposition gives us a way to \emph{disprove} well-posedness in coarse topologies, simply by establishing that at least one of the
operators $A_n$ is discontinuous.  This tends to be the case if $A_n$ contains a significant ``high-to-low frequency cascade'',
and we shall exploit this (in the $n=2$ case) to prove Theorem \ref{ill-p}.

\section{Reduction to function spaces}\label{lwp-sec}

We can now reduce Theorem \ref{lwp} and Theorem \ref{ill-p} to the construction of a certain pair $\S^s(\R \times \R)$ and $\N^s(\R \times \R)$
of function spaces for each $s \in \R$, the verification of certain estimates for these spaces, and the verification of a certain bad behavior
of the quadratic operator $A_2$. The precise statements are as follows. Define a \emph{bump function} to be
a smooth compactly supported function $t \mapsto \eta(t)$ of the time variable $t \in \R$.

\begin{proposition}[Function spaces]\label{func-sec}  For any $s \in \R$ there exists a Banach space $\S^s(\R \times \R)$ (the ``solution space'' at regularity $H^s_x(\R)$) and a Banach space $\N^s(\R \times \R)$ (the ``nonlinearity space'' at regularity $H^s_x(\R)$), with the following properties:
\begin{itemize}
\item[(i)] (Density) The Schwartz functions on $\R \times \R$ are dense in $\S^s(\R \times \R)$ and in $\N^s(\R \times \R)$.
\item[(ii)] (Nesting)  If $s \leq s'$ and $u \in \S^{s'}(\R \times \R)$, then
$$ \| u \|_{\S^s(\R \times \R)} \leq \| u \|_{\S^{s'}(\R \times \R)}.$$
Similarly, if $F \in \N^{s'}(\R \times \R)$ then
$$ \| F \|_{\N^s(\R \times \R)} \leq \| F \|_{\N^{s'}(\R \times \R)}.$$
\item[(iii)] (Energy estimate) If $u \in \S^s(\R \times \R)$, then\footnote{We use $C_s$ to denote a positive constant - which can vary from line to line - that can depend on $s$.  Similarly if we subscript $C$ by other parameters.}
$$ \| u \|_{C^0_t H^s_x(\R \times \R)} \leq C_s \| u \|_{\S^s(\R \times \R)}.$$
\item[(iv)] (Homogeneous estimate)  If $u_0 \in H^s_x(\R)$, $u(t) = \exp(i t \partial_{xx}) u_0$, 
and $\eta(t)$ is a bump function, then $\eta u \in \S^s(\R \times \R)$ and
$$ \| \eta u \|_{\S^s(\R \times \R)} \leq C_{\eta,s} \| u_0 \|_{H^s_x(\R)}.$$
\item[(v)] (Dual homogeneous estimate)  If $F \in \N(\R \times \R)$, and $\eta(t)$ is a bump function, then
$$ \| \int_\R \sgn(s) \eta(s) \exp(-is\partial_{xx}) F(s)\ ds \|_{H^s_x(\R)} \leq C_{\eta,s} \| F \|_{\N^s(\R \times \R)}.$$
\item[(vi)] (Inhomogeneous estimate) If $F \in \N^s(\R \times \R)$, and $\eta(t)$ is a bump function, then
$$ \| \int_\R \sgn(t-s) \eta(t-s) \exp(-i(t-s)\partial_{xx}) F(s)\ ds \|_{\S^s(\R \times \R)} \leq C_{\eta,s} \| F \|_{\N^s(\R \times \R)}.$$
\item[(vii)] (Nonlinear estimate)  If $u,v \in \S^s(\R \times \R)$ for some $s \geq -1$, then
$$ \| uv \|_{\N^s(\R \times \R)} \leq 
C_{s}
( \| u \|_{\S^s(\R \times \R)} \| v \|_{\S^{-1}(\R \times \R)} + \| u \|_{\S^{-1}(\R \times \R)} \| v \|_{\S^s(\R \times \R)} ).$$
Here we define the non-linear operation $(u,v) \mapsto uv$ first for Schwartz functions, and then extend to the general case by density.
\end{itemize}
\end{proposition}

We shall prove this proposition in later sections.  Assuming it for now, we see from (iv)-(vii) and
\eqref{Ldef}, \eqref{Ndef} that the equation \eqref{uln-2} is quantitatively well-posed in the spaces $H^{-1}_x(\R)$, $\S^{-1}(\R \times \R)$,
and also one has persistence of regularity (in the sense of Theorem \ref{persist} for the spaces $H^s_x(\R)$, $\S^s(\R \times \R)$
for any $s \geq -1$.  Combined with (ii) and (iii), we thus establish Theorem \ref{lwp} (using the reductions in Section \ref{reduced-sec}).

The derivation of Theorem \ref{ill-p} from this Proposition is almost as immediate.

\begin{proof}[Proof of Theorem \ref{ill-p}]
Fix $s < -1$ and $s'$; we may rescale $T$ to equal $1$.  Suppose for contradiction
that solution map $f \mapsto u[f]$ is continuous on $B_r$ (with the $H^s_x(\R)$ topology) to $C^0_t H^{-1}_x([0,1] \times \R)$ 
(with the $C^0_t H^{s'}_x$ topology).  Applying Proposition \ref{power}, we conclude that the quadratic operator
$$ A_2: f \mapsto N_2(Lf, Lf)$$ 
(restricted of course to $[0,1] \times \R$) is continuous from $B_r$ (with the $H^s_x(\R)$ topology) 
to $C^0_t H^{s'}_x([0,1] \times \R)$.  In particular, this implies that
$$ \sup_{0 \leq t\leq 1}\| A_2(f_N)(y) \|_{H^{s'}(\R)} \to 0$$
whenever $f_N$ is a sequence of functions in $B_r$ which goes to zero in $H^s_x(\R)$ norm.
The left-hand side can be expanded by \eqref{Ldef}, \eqref{Ndef} as
$$ \sup_{0 \leq t \leq 1} 
\| \int_0^t \exp(i(t-t')\partial_{xx})((\exp(it' \partial_{xx}) f_N)^2)\ dt' \|_{H^{s'}(\R)}$$
which after taking Fourier transforms becomes
$$\sup_{0 \leq t \leq 1} 
\| \langle \xi \rangle^{s'}
\int_0^t \int_\R \exp(-i(t-t')\xi^2) \exp(it' (\xi_1^2 + (\xi-\xi_1)^2) \hat f_N(\xi_1) \hat f_N(\xi-\xi_1)\ d\xi_1 dt'
\|_{L^2_\xi(\R)}.$$
Now let $N > 100$ be a large parameter, and set 
$$ \hat f_N := r N 1_{[-10,10]}(|\xi|-N) / 1000.$$
Then $f_N \in B_r$, and $\|f_N\|_{H^s_x(\R)} \to 0$ as $N \to \infty$.  Thus
\begin{equation}\label{sin}
\sup_{0 \leq t \leq 1} \| \langle \xi \rangle^{s'}
\int_0^t \int_\R \exp(-i(t-t')\xi^2) \exp(it' (\xi_1^2 + (\xi-\xi_1)^2) \hat f_N(\xi_1) \hat f_N(\xi-\xi_1)\ d\xi_1 dt'
\|_{L^2_\xi(\R)} \to 0
\end{equation}
as $N \to \infty$.
Now set $t := 1/100N^2$ and localize to the region where $-1 \leq \xi \leq 1$.  One can verify that
$$ \Re( \exp(-i(t-t')\xi^2) \exp(it' (\xi_1^2 + (\xi-\xi_1)^2) ) > 1/2$$
whenever $0 \leq t' \leq t$ and $\xi_1$ lives in the support of $f$, hence we obtain
$$ 
\| \langle \xi \rangle^{s'}
\int_0^t \int_\R \exp(-i(t-t')\xi^2) \exp(it' (\xi_1^2 + (\xi-\xi_1)^2) \hat f_N(\xi_1) \hat f_N(\xi-\xi_1)\ d\xi_1 dt'
\|_{L^2_\xi(\R)} \geq c r^2$$
for some $c > 0$.  But this contradicts \eqref{sin}.
This proves Theorem \ref{ill-p}.
\end{proof}

It is instructive to compute the spacetime Fourier transform of $A_2(f)$.  A computation shows that it has a significant component near the time-frequency axis $\xi = 0$, indeed it has magnitude comparable to $N^{-2s-1}$ on
a rectangle $\{ \xi = O(1), \tau = 2N^2 + O( N ) \}$. This is already enough to cause it to leave the $X^{s,b}$ space whenever $s < b - \frac{5}{4}$, which explains why the $X^{s,b}$ method ceases to work well when $s < -3/4$.  However, this will be fixed by replacing this space with an $L^1_\tau$-based space near the time axis.

\section{Taking the spacetime Fourier transform}

To complete the proof of Theorem \ref{lwp} and Theorem \ref{ill-p}, we need to build spaces $\S^s(\R \times \R)$ and $\N^s(\R \times \R)$
which obey the properties in Proposition \ref{func-sec}.  To abbreviate the notation we shall now omit the domain $\R \times \R$ from
these spaces.  For $s > 1/2$ one could use the energy spaces $\S^s := C^0_t H^s_x$,
$\N^s := L^1_t H^s_x$, while for $s \geq 0$ one could use Strichartz spaces such as $\S^s := C_0^t H^s_x \cap L^4_t L^\infty_x$,
$\N^s := L^1_t H^s_x$ (other choices are available; see \cite{tsutsumi}, \cite{cwi}).  For $s > -3/4$, it was shown in \cite{kpv}
that one could use the spaces $\S^s := X^{s,b}$ and $\N^s := X^{s,b-1}$ for any $b > 1/2$.

In order to construct spaces which work all the way down to $s \geq -1$, we have to modify the $X^{s,b}$ spaces somewhat.  The precise
modification is somewhat complicated, so for now we shall continue to work abstractly to avoid being bogged down in details.  We shall require
a space $W$, constructed by the following proposition.  Call a function on $\R \times \R$ \emph{reasonable} if it lies in $L^\infty_t L^\infty_x(\R \times \R)$ and has compact support.  For any $s,b\in \R$, we define $\hat X^{s,b}$ to be the closure of the reasonable functions via the norm
$$ \|f\|_{\hat X^{s,b}} := \| \langle \xi \rangle^s \langle \tau - \xi^2 \rangle^b f\|_{L^2_\tau L^2_\xi}.$$
These are the Fourier transforms of the usual $X^{s,b}$ spaces.

\begin{proposition}[Construction of main space]\label{construct}  There exists a Banach space $W$, which is the closure of the reasonable
functions in $\R \times \R$ by some norm $\| \|_W$, with the following properties for all reasonable
$f,g$:
\begin{itemize}
\item (Monotonicity) If $|f| \leq |g|$ pointwise, then $\|f\|_W \leq \|g\|_W$.  In particular, we have $\|f\|_W = \||f|\|_W$ (so the $W$ norm depends only on the magnitude of the function, and not the phase).
\item ($H^{-1}$ Energy estimate)  We have
\begin{equation}\label{W-energy}
\| \langle \xi \rangle^{-1} f \|_{L^2_\xi L^1_\tau} \leq C \|f\|_W
\end{equation}
\item (Homogeneous $H^{-1}$ solution estimate)  We have
\begin{equation}\label{W-homogeneous}
\| f \|_W \leq C \| f \|_{\hat X^{-1,100}}.
\end{equation}
\item (Bilinear estimate) We have
\begin{equation}\label{W-bilinear}
\| \langle \tau - \xi^2 \rangle^{-1} f * g \|_W
\leq \| f \|_W \|g\|_W,
\end{equation}
where of course $f*g$ denotes spacetime convolution
$$ f*g(\tau,\xi):= \int_\R \int_\R f(\tau_1,\xi_1) g(\tau_2,\xi_2)\ d\tau_1 d\xi_1$$
using the convention 
\begin{equation}\label{tau-convention}
(\tau_1,\xi_1) + (\tau_2,\xi_2) = (\tau,\xi).
\end{equation}  
\end{itemize}
\end{proposition}

The space $W$ has been designed with the scaling of $H^{-1}_x$, as this is the most important regularity in our argument.  
A good candidate to keep in mind for $W$ is the space $\hat X^{-1,b}$ for some $b > 1/2$; this turns out to only obey the first
three properties required (and enough of the fourth property that one can establish local existence in $s > -3/4$ rather than $s \geq -1$);
our final version of $W$ shall be a modification of $\hat X^{-1,b}$. 

We prove this proposition in the remainder of the paper.  For now, let assume it, and use it to prove Proposition \ref{func-sec}.  We will
take $\S^s$ and $\N^s$ to be the closure of the Schwartz functions under the norms
\begin{equation}\label{snef}
\| u \|_{\S^s} := \| \langle \xi \rangle^{s+1} \tilde u \|_W; \quad \| F \|_{\N^s} := \left\| \frac{\langle \xi \rangle^{s+1}}{\langle \tau - \xi^2 \rangle} \tilde F \right\|_W
\end{equation}
where $\tilde u(\tau,\xi)$ denotes the spacetime Fourier transform.

By construction, the density and nesting properties of $\S^s$ and $\N^s$ required for Proposition \ref{func-sec} are immediate.
To prove the energy estimate (iii), it suffices by the usual limiting arguments to show that
$$ \sup_t \| u(t) \|_{H^s_x} \leq C_s \| u \|_{\S^s}$$
when $u$ is Schwartz.
From Fourier inversion and the triangle inequality we have
$$ \| u(t) \|_{H^s_x} = \left\| \int_\R \langle \xi \rangle^s e^{it\tau} \tilde u(\tau,\xi)\ d\tau \right\|_{L^2_\xi}
\leq \| \langle \xi \rangle^s \tilde u \|_{L^2_\xi L^1_\tau}$$
and the claim now follows from \eqref{W-energy}.

To prove the homogeneous estimate (iv), observe that if $u(t) = \exp(it\partial_{xx}) u_0$, then
$$ \widetilde{\eta u}(\tau,\xi) = \hat \eta(\tau - \xi^2) \hat u_0(\xi)$$
and in particular from the rapid decrease of $\hat \eta$ we see that
$$ \| \eta u \|_{\hat X^{s,100}} \leq C_\eta \| u_0 \|_{H^s_x}$$
(say).  Thus the claim follows from \eqref{W-homogeneous}.

To prove the dual homogeneous estimate (v), we apply the Fourier transform and Parseval's identity in space to write
$$ \| \int_\R \sgn(s) \eta(s) \exp(-is\partial_{xx}) F(s)\ ds \|_{H^s_x}
= \sqrt{2\pi} \| \int_\R \langle \xi \rangle^s \tilde F(\tau,\xi) \widehat{\sgn \eta}(\tau-\xi^2)\ d\tau \|_{L^2_\xi}.$$
A simple integration by parts shows that 
\begin{equation}\label{wet}
|\widehat{\sgn \eta}(\tau)| \leq C_\eta \langle \tau \rangle^{-1}.
\end{equation}
Inserting this bound and then using \eqref{W-energy}, \eqref{snef} we obtain
$$ \| \int_\R \sgn(s) \eta(s) \exp(-is\partial_{xx}) F(s)\ ds \|_{H^s_x}
\leq C_\eta \| \langle \tau - \xi^2\rangle^{-1}\langle \xi \rangle^s |\tilde F| \|_W \leq C_\eta \| F \|_{\N^s}$$
as desired.

To establish the inhomogeneous estimate, observe that the spacetime Fourier transform of
$\int_\R \sgn(t-s) \eta(t-s) \exp(-i(t-s)\partial_{xx}) F(s)\ ds$ at $(\tau,\xi)$ is simply
$\widehat{\sgn \eta}(\tau - \xi) \tilde F(\tau,\xi)$.  Using \eqref{wet} and \eqref{snef} we obtain
$$ \| \int_\R \sgn(t-s) \eta(t-s) \exp(-i(t-s)\partial_{xx}) F(s)\ ds \|_{\S^s} \leq C_{\eta} \| F \|_{\N^s}$$
as desired.

Finally, we consider the non-linear estimate.  Let us write $s = -1 + \delta$.  Taking Fourier transforms, we have
$$ \| uv \|_{\N^{-1+\delta}} = 
\left\| \int_\R \int_\R \frac{\langle \xi \rangle^\delta}{\langle \tau - \xi^2 \rangle}
\tilde u(\tau_1,\xi_1) \tilde v(\tau_2,\xi_2)\ d\xi_1 d\tau_1
\right\|_W$$
where $(\tau_2,\xi_2)$ is defined by the convention \eqref{tau-convention}.
We then estimate
$$ \langle \xi \rangle^\delta \leq C_s( \langle \xi_1 \rangle^\delta + \langle \xi_2 \rangle^\delta )$$
and so by symmetry it would suffice to show that
$$
\| \int_\R \int_\R \frac{\langle \xi_1 \rangle^\delta}{\langle \tau - \xi^2 \rangle}
|\tilde u(\tau_1,\xi_1)| |\tilde v(\tau_2,\xi_2)|\ d\xi_1 d\tau_1 \|_W
\leq \| u \|_{\S^{-1+\delta}} \| v \|_{\S^{-1}}.$$
Writing $|\tilde u(\tau,\xi)| = \langle \xi \rangle^{-\delta} f(\tau,\xi)$ and $|\tilde v(\tau,\xi)| = g(\tau,\xi)$, 
the claim then follows from \eqref{W-bilinear}.

It remains to prove Proposition \ref{construct}.  For now, let us make an important consequence of
the convention \eqref{tau-convention} which is essential to the argument.  Observe that 
\eqref{tau-convention} implies the algebraic identity
$$ (\tau - \xi^2) = (\tau_1 - \xi_1^2) + (\tau_2 - \xi_2^2) - 2 \xi_1 \xi_2$$
and so by the triangle inequality we have the fundamental \emph{resonance estimate}
\begin{equation}\label{resonance}
\max( \langle \tau - \xi^2 \rangle, \langle \tau_1 - \xi_1^2 \rangle, \langle \tau_2 - \xi_2^2 \rangle)
\geq 2^{-5} \langle \xi_1 \xi_2 \rangle.
\end{equation}
(The constant $2^{-5}$ is very conservative, but the exact value of it is not important here.)
Thus if both input frequencies $\xi_1$ and $\xi_2$ are large, then it is not possible for all three of
$(\tau,\xi^2)$, $(\tau_1,\xi_1^2)$, and $(\tau_2,\xi_2^2)$ to lie close to the parabola.

\section{Description of the space $W$}

We are now ready to construct $W$ and establish all the desired properties except for the bilinear
estimate \eqref{W-bilinear}, which will be straightforward but require some effort. As mentioned earlier,
the space $\hat X^{-1,b}$ is the model candidate for $W$.  However we need to make three modifications to this
space in order for it to be viable for us all the way down to the endpoint $s=-1$.

It will not be surprising that the geometry of the parabola $\tau = \xi^2$ plays a crucial role.  We shall rely in particular
on localizations to the spatial annuli
$$ A_j := \{ (\tau,\xi) \in \R \times \R: 2^j \leq \langle \xi \rangle < 2^{j+1} \}$$
for $j \geq 0$, as well as the parabolic neighborhoods
$$ B_d := \{ (\tau,\xi) \in \R \times \R: 2^d \leq \langle \tau - \xi^2 \rangle < 2^{d+1} \}$$
for $d \geq 0$.  Thus the sets $A_j \cap B_d$ for $j,d \geq 0$ partition frequency space, and we have\footnote{All sums and unions involving $j$ and $d$ shall be over the non-negative integers unless otherwise mentioned.}
\begin{equation}\label{wt}
 \| f\|_{\hat X^{s,b}} \approx 
(\sum_{j} \sum_{d} 2^{2sj} 2^{2bd} \| f \|_{L^2_\xi L^2_\tau(A_j \cap B_d)}^2)^{1/2}.
\end{equation}
We also use the variant sets
$$ A_{\leq j} := \bigcup_{j' \leq j} A_{j'}; \quad B_{\leq d} := \bigcup_{d' \leq d} B_{d'}$$
and similarly define $A_{\geq j}$, $A_{>j}$, $B_{\geq d}$, $B_{>d}$, etc.

A natural candidate for the space $W$ is then the Besov endpoint $\hat X^{-1,1/2,1}$ of \eqref{wt}, defined by
\begin{equation}\label{fx}
 \| f \|_{\hat X^{-1,1/2,1}} :=  (\sum_{j} 2^{-j} (\sum_{d} 2^{d/2} \| f \|_{L^2_\xi L^2_\tau(A_j \cap B_d)})^2)^{1/2}.
\end{equation}
This type of space has appeared previously in endpoint theory (see for instance \cite{tataru:wave1}); we shall
need this Besov refinement in order to handle the $s=-1$ endpoint without encountering logarithmic divergences (in particular, to handle
the ``parallel interaction'' case when the nonlinearity interacts two components of the solution $u$ with the same high frequency).
The relationship between this space and the $\hat X^{s,b}$ spaces are provided by the following easy lemma.

\begin{lemma}  Suppose $f$ supported on $B_{\geq d}$ for some $d \geq 0$ (this condition is vacuous when $d=0$).  Then we have
\begin{equation}\label{x-trivial-left}
 \| f \|_{\hat X^{-1,b}} \leq C_b 2^{-(1/2-b)d} \|f\|_{\hat X^{-1,1/2,1}}
\end{equation}
whenever $b < 1/2$, and
\begin{equation}\label{x-trivial-right}
 \| f \|_{\hat X^{-1,1/2,1}} \leq C_b 2^{-(b-1/2)d} \|f\|_{\hat X^{-1,b}}
\end{equation}
whenever $b > 1/2$.
\end{lemma}

\begin{proof}  We may easily restrict $f$ to a single annulus $A_j$, since the general case then follows by square-summing.  The claim
then follows by decomposing further into $A_j\cap B_{d'}$ for $d' \geq d$ and using Cauchy-Schwarz.
\end{proof}

By using $\hat X^{-1,1/2,1}$ instead of $\hat X^{-1,b}$, we will be able to handle parallel interactions.  However, as essentially 
observed in \cite{mur}, this Besov refinement is not sufficient by itself even to handle the endpoint $s=-3/4$, because
of a divergence at the time axis $\tau=0$.  To handle these divergences we need a somewhat different norm $\| \|_{Y}$, defined as
\begin{equation}\label{fy}
\|f\|_{Y} := \| \langle \xi \rangle^{-1} f \|_{L^2_\xi L^1_\tau} + \| f \|_{L^2_\xi L^2_\tau} 
\end{equation}
and then form the sum space $Z := \hat X^{-1,1/2,1}+Y$ in the usual fashion,
$$ \| f \|_{Z} := \inf \{ \| f_1 \|_{\hat X^{-1,1/2,1}} + \|f_2\|_{Y}: f_1 \in \hat X^{-1,1/2,1}; f_2 \in Y; f=f_1+f_2 \}.$$
It is easy to verify that this is a Banach space, with the Schwartz functions being dense.  Clearly we have
$$ \|f\|_{Z} \leq \|f\|_{\hat X^{-1,1/2,1}}, \|f\|_Y ,$$
and conversely, to prove any linear estimate of the form $\|Tf\|_Z \leq M \|f\|_{Z}$ where $Z$ is a Banach space and $M > 0$, it suffices to prove the separate estimates $\|Tf\|_Z \leq M\|f\|_{\hat X^{-1,1/2,1}}$, $\|Tf\|_Z \leq M \|f\|_Y$.  For instance, we have
the following basic estimates.

\begin{proposition}\label{space-est}  For any reasonable $f$, we have
\begin{align}
\|\langle \xi \rangle^{-1} f\|_{L^2_\xi L^1_\tau} &\leq C \|f\|_{Z}\label{f21}
\end{align}
If furthermore $f$ is supported on $A_j \cap B_{\geq d}$ for some $j, d \geq 0$, then we have
\begin{align}
\|f\|_{L^2_\xi L^2_\tau} &\leq C (1 + 2^{j} 2^{-d/2}) \|f\|_{Z}.\label{f22}\\
\|f\|_{L^1_\xi L^1_\tau} &\leq C 2^{3j/2} \|f\|_{Z}\label{f11} \\
\|f\|_{L^1_\xi L^2_\tau} &\leq C (2^{j/2} + 2^{3j/2} 2^{-d/2}) \|f\|_{Z}\label{f12}.
\end{align}
\end{proposition}

\begin{proof}  Observe that \eqref{f21}, \eqref{f22} follow immediately from \eqref{fy} if the right-hand sides were replaced by $C \|f\|_Y$, so it
suffices to establish these estimates with the right-hand side of $C \|f\|_{\hat X^{-1,1/2,1}}$.  For \eqref{f22}, this follows from
\eqref{x-trivial-left}.  As for \eqref{f11}, we may reduce to a single annulus $A_j$ (after square-summing in $j$) and reduce
to showing that
$$ \| f\|_{L^2_\xi L^1_\tau(A_j)} \leq C \sum_{d} 2^{d/2} \| f \|_{L^2_\xi L^2_\tau(A_j \cap B_d)}.$$
But this follows from the triangle inequality and H\"older's inequality, since for each fixed $\xi$ and fixed $B_d$, the $\tau$ variable varies
over a set of measure $O( 2^d )$.

Finally, \eqref{f11}, \eqref{f12} follow respectively from \eqref{f21}, \eqref{f22} and H\"older in the $\xi$ variable (which varies over a set of measure $O(2^j)$).
\end{proof}

Observe that $\hat X^{-1,1/2,1}$ and $Y$ are both monotone in the sense of Proposition \ref{construct}.  Also
the two spaces $\hat X^{-1,1/2,1}$ and $Y$ paste together nicely along the fuzzy boundary $\langle \tau - \xi^2 \rangle \approx \langle \xi \rangle^2$.
A formalization of this heuristic is as follows.

\begin{lemma}[Pasting lemma]\label{paste}  Let $f$ be a reasonable function.  If $f$ is supported on $\bigcup_{j} A_j \cap B_{\geq 2j - 100}$, then
\begin{equation}\label{paste-a}
\|f\|_Y \leq C \|f\|_{Z}.
\end{equation}
Conversely, if $f$ is supported on $\bigcup_{j} A_j \cap B_{\leq 2j + 100}$, then
\begin{equation}\label{paste-b}
 \|f\|_{\hat X^{-1,1/2,1}} \leq C \|f\|_{Z}.
 \end{equation}
\end{lemma}

\begin{proof}  Let us first establish \eqref{paste-a}.  It clearly suffices to show that
$$ \|f\|_Y \leq C \|f\|_{\hat X^{-1,1/2,1}}$$
on this domain.  Partitioning dyadically into the $A_j$, and then square-summing in $j$, it suffices to show that
$$ \|f 1_{A_j}\|_Y \leq C \|f 1_{A_j}\|_{\hat X^{-1,1/2,1}}$$
for each $j$.  From \eqref{x-trivial-left} we already have
$$ \|f 1_{A_j}\|_{L^2_\xi L^2_\tau} \leq C \|f 1_{A_j}\|_{\hat X^{-1,1/2,1}}$$
while from \eqref{f21} we have
$$ \|\langle \xi \rangle^{-1} f 1_{A_j}\|_{L^2_\xi L^1_\tau} \leq C \|f 1_{A_j}\|_{\hat X^{-1,1/2,1}}$$
and the claim follows.

Now we establish \eqref{paste-b}.  By arguing as before it suffices to show that
$$  \|f 1_{A_j}\|_{\hat X^{-1,1/2,1}} \leq C \|f 1_{A_j}\|_{Y}.$$
But we have
\begin{align*}
 \|f 1_{A_j}\|_{\hat X^{-1,1/2,1}} &\leq C \sum_{d \leq 2j+100} 2^{-j} 2^{d/2} \| f \|_{L^2_\xi L^2_\tau(A_j \cap B_d)} \\
&\leq C \| f \|_{L^2_\xi L^2_\tau(A_j)}\\
&\leq C \|f 1_{A_j} \|_Y
\end{align*}
as desired.
\end{proof}

The space $Z$ is a candidate for $W$, as it is able to cope with two of the dangerous quadratic interactions in 
the equation (namely the parallel interactions, and the interactions which output near the time axis).  However, there is a third type 
of interaction which could cause trouble, when a solution component near the parabola $\{ \tau = \xi^2 \}$ interacts with a component near the reflected parabola $\{ \tau = -\xi^2 \}$ to create a large contribution
near the frequency origin.  The use of the space $Z$ does not prohibit either component from occuring.  However, as mentioned 
in the introduction, the solution should stay in the upper half-plane $\tau > 0$.  To exploit this we shall introduce a weight 
\begin{equation}\label{weight-def}
w(\tau,\xi) := \max(1,-\tau)^{10}
\end{equation}
to localize to the upper half-plane, and define $W$ to be the space
\begin{equation}\label{W-def}
 \|f\|_W := \| w f \|_{Z}
\end{equation}
as discussed in the previous section.

The monotonicity of $W$ is clear.  
The claim \eqref{W-energy} follows immediately from \eqref{f21} (since $w \geq 1$), while \eqref{W-homogeneous} follows from
\eqref{x-trivial-right}:
$$ \|f\|_W \leq \|w f\|_{\hat X^{-1,1/2,1}} \leq C\| w f\|_{\hat X^{-1,90}} \leq C\|f\|_{\hat X^{-1,100}}$$
where we use the crude estimate $w(\tau,\xi) \leq C \langle \tau - \xi^2 \rangle^{10}$.

It remains to show \eqref{W-bilinear}.  Applying \eqref{W-def} and monotonicity, we reduce to showing that
\begin{equation}\label{main-bilinear}
 \| \frac{w}{\langle \tau - \xi^2 \rangle} (\frac{f}{w} * \frac{g}{w}) \|_{Z}
\leq C \|f\|_{Z} \|g\|_{Z}.
\end{equation}  
for all non-negative reasonable $f,g$.

Fix $f,g$.
Observe from \eqref{tau-convention} that
$$ w(\tau,\xi) \leq C w(\tau_1,\xi_1) w(\tau_2,\xi_2);$$
this basically reflects the fact that in order for $\tau = \tau_1+\tau_2$ to be negative, at least one of $\tau_1,\tau_2$ has to be even more negative.
This gives us the very handy pointwise estimate
\begin{equation}\label{k-point}
\frac{w}{\langle \tau - \xi^2 \rangle} (\frac{f}{w} * \frac{g}{w}) \leq \frac{C}{\langle \tau - \xi^2 \rangle} (f*g),
\end{equation}
which we shall rely upon in most cases (except for one special high-high interaction where we must utilize the localizing weight $w$ more carefully).  
As one example of this, we present a relatively simple case of \eqref{main-bilinear}:

\begin{lemma}\label{yy-bilinear}  For any non-negative reasonable $f,g$, we have
$$  \| \frac{w}{\langle \tau - \xi^2 \rangle} (\frac{f}{w} * \frac{g}{w}) \|_{Z}
\leq C \|f\|_{Y} \|g\|_{Y}.$$
\end{lemma}

\begin{proof}  Let us first restrict the left-hand side to the region $|\tau| \leq 10 \langle \xi^2 \rangle$.  From Young's inequality we have
$$ \|f*g\|_{L^\infty_\tau L^\infty_\xi} \leq \|f\|_{L^2_\tau L^2_\xi} \|g\|_{L^2_\tau L^2_\xi} \leq \|f\|_Y \|g\|_Y$$
so by \eqref{k-point} and monotonicity it suffices to show that
$$ \| \frac{1}{\langle \tau - \xi^2 \rangle} 1_{|\tau| \leq 10 \langle \xi^2\rangle}
\|_{\hat X^{-1,1/2,1}} \leq C.$$
But this is easily verified.

Now we turn to the region where $|\tau| > 10 \langle \xi^2 \rangle$.  In this regime $\langle \tau - \xi^2 \rangle$ is comparable to $\langle \tau\rangle$, and so it will suffice to show that
$$ \left( \sum_j 2^{-j} [\sum_d 2^{d/2} \| \frac{w}{2^d} (\frac{f}{w} * \frac{g}{w}) \|_{L^2_\xi L^2_\tau(A_j \cap C_d)}]^2 \right)^{1/2} \leq C \|f\|_{L^2_\xi L^2_\tau} \|g\|_{L^2_\xi L^2_\tau}$$
where $C_d := \{ (\tau,\xi): 2^d < \langle \tau\rangle < 2^{d+1} \}$.

By H\"older's inequality, it will be enough to show that
$$ \sum_d 2^{-d/2} \| w (\frac{f}{w} * \frac{g}{w}) \|_{L^\infty_\xi L^2_\tau(C_d)} \leq C \|f\|_{L^2_\xi L^2_\tau} \|g\|_{L^2_\xi L^2_\tau}.$$

By the triangle inequality, we can bound the left-hand side by
$$ \sum_d \sum_{d_1} \sum_{d_2} 2^{-d/2} \| w (\frac{f 1_{C_{d_1}}}{w} * \frac{g 1_{C_{d_2}}}{w}) \|_{L^\infty_\xi L^2_\tau(C_d)}.$$
We may assume $\max(d_1,d_2) \geq d-10$, since the summand vanishes otherwise.

Suppose first that $d-10 \leq d_1 \leq d+10$ and $d_2 \leq d+10$.  Then we use \eqref{k-point} and Young's inequality to bound
\begin{align*}
\| w (\frac{f 1_{C_{d_1}}}{w} * \frac{g 1_{C_{d_2}}}{w}) \|_{L^\infty_\xi L^2_\tau(C_d)} 
&\leq C \| f 1_{C_{d_1}} * g 1_{C_{d_2}} \|_{L^\infty_\xi L^2_\tau(C_d)} \\
&\leq C \| f \|_{L^2_\xi L^2_\tau(C_{d_1})} \| g \|_{L^2_\xi L^1_\tau(C_{d_2})} \\
&\leq C 2^{d_2/2} \| f \|_{L^2_\xi L^2_\tau(C_{d_1})} \| g \|_{L^2_\xi L^2_\tau(C_{d_2})} 
\end{align*}
and the contribution of this case can be disposed of by Schur's test.  A similar argument deals with the case when $d-10 \leq d_2 \leq d+10$ and $d_1 \leq d+10$.  Thus it only remains to consider the case $\max(d_1,d_2) > d+10$, in which case $|d_1-d_2| \leq 2$.  But in this regime one can strengthen \eqref{k-point} to
$$ w (\frac{f 1_{C_{d_1}}}{w} * \frac{g 1_{C_{d_2}}}{w}) \leq C 2^{-10 d} f 1_{C_{d_1}} * g 1_{C_{d_2}}.$$
Using this and Young's inequality, we have
\begin{align*}
\| w (\frac{f 1_{C_{d_1}}}{w} * \frac{g 1_{C_{d_2}}}{w})  \|{L^\infty_\xi L^2_\tau(C_d)} 
&\leq C 2^{d/2} 2^{-10d}
\| f 1_{C_{d_1}} * g 1_{C_{d_2}} \|{L^\infty_\xi L^\infty_\tau(C_d)} \\
&\leq  C 2^{d/2} 2^{-10d} \| f \|_{L^2_\xi L^2_\tau(C_{d_1})} \| g \|_{L^2_\xi L^2_\tau(C_{d_2})}
\end{align*}
and the contribution of this case can again be dealt with by Schur's test.
\end{proof}

From this lemma, Lemma \ref{paste}, and monotonicity, we thus see that to prove \eqref{main-bilinear} we may thus assume
that at least one of non-negative reasonable $f,g$ lies near the parabola, or more precisely
\begin{equation}\label{assumption}
\hbox{At least one of } f, g \hbox{ is supported in } \bigcup_{j} A_j \cap B_{< 2j - 100}.
\end{equation}

The next step is dyadic decomposition.  
Observe the localization property
\begin{equation}\label{localization}
\|f\|_{Z} \sim (\sum_j\|1_{A_j} f \|_{Z}^2)^{1/2},
\end{equation}
which follows from $L^2_\xi$ nature of both $\hat X^{-1,1/2,1}$ and $Y$.  We therefore
split $f = \sum_{j_1} f_{j_1}$ and
$g = \sum_{j_2} g_{j_2}$, where $f_{j_1}$ and $g_{j_2}$ are the restrictions of $f$, $g$ to $A_{j_1}$, $A_{j_2}$ respectively.
Thus
\begin{align*}
 \| \frac{w}{\langle \tau - \xi^2 \rangle} (\frac{f}{w} * \frac{g}{w}) \|_{Z}
 &\leq C
  (\sum_j \| 1_{A_j} \frac{w}{\langle \tau - \xi^2 \rangle} (\frac{f}{w} * \frac{g}{w}) \|_{Z}^2)^{1/2}\\
 &= C
   (\sum_j \| \sum_{j_1,j_2} 1_{A_j} \frac{w}{\langle \tau - \xi^2 \rangle} (\frac{f_{j_1}}{w} * \frac{g_{j_2}}{w}) \|_{Z}^2)^{1/2}.
\end{align*}
In order for the inner summand to be non-zero, it must be possible to find $(\tau,\xi) \in A_j$, $(\tau_1,\xi_1) \in A_{j_1}$, $(\tau_2, \xi_2) \in A_{j_2}$ obeying \eqref{tau-convention}.  This forces one of the following (overlapping) cases to hold:

\begin{itemize}
\item (High-low interaction) $|j-j_1| \leq 10$ (which implies $j_2 \leq j+11$);
\item (Low-high interaction) $|j-j_2| \leq 10$ (which implies $j_1 \leq j+11$);
\item (High-high interaction) $j < j_1-10, j_2-10$ (which implies $|j_1-j_2| \leq 1$).
\end{itemize}

The former two cases are symmetric.  Thus to prove \eqref{main-bilinear} it suffices (again using \eqref{localization}) to verify the high-low estimate
\begin{align*}
&(\sum_j \| \sum_{|j_1-j| \leq 10} \sum_{j_2 \leq j+11} 
1_{A_j} \frac{w}{\langle \tau - \xi^2 \rangle} (\frac{f_{j_1}}{w} * \frac{g_{j_2}}{w}) \|_{Z}^2)^{1/2}\\
&\quad\leq C (\sum_{j_1} \|f_{j_1}\|_{Z}^2)^{1/2} (\sum_{j_2} \|g_{j_2}\|_{Z}^2)^{1/2}
\end{align*}
and the high-high estimate
\begin{equation}\label{hhest}
\begin{split}
&(\sum_j \| \sum_{j_1,j_2 > j+10: |j_1-j_2| \leq 1} 1_{A_j} \frac{w}{\langle \tau - \xi^2 \rangle} (\frac{f_{j_1}}{w} * \frac{g_{j_2}}{w}) \|_{Z}^2)^{1/2} \\
&\quad \leq  C (\sum_{j_1} \|f_{j_1}\|_{Z}^2)^{1/2} (\sum_{j_2} \|g_{j_2}\|_{Z}^2)^{1/2}.
\end{split}
\end{equation}
Consider the high-low estimate first.  We use \eqref{k-point} to drop the weights $w$.  Since for any $j$ there are only $O(1)$ values of $j_1$ which contribute, we can use Schur's test and reduce to showing
$$
\| \sum_{j_2 \leq j+11} 
\frac{1_{A_j}}{\langle \tau - \xi^2 \rangle} (f_{j_1} * g_{j_2}) \|_{Z}
\leq C \|f_{j_1}\|_{Z} (\sum_{j_2} \|g_{j_2}\|_{Z}^2)^{1/2}$$
whenever $|j_1-j| \leq 10$.  By the triangle inequality, and estimating the $Z$ norm by
the $\hat X^{-1,1/2,1}$ norm, it thus suffices to establish
\begin{equation}\label{high-low-targ}
\| \frac{1_{A_j}}{\langle \tau - \xi^2 \rangle} (f_{j_1} * g_{j_2}) \|_{\hat X^{-1,1/2,1}}
\leq C (2^{-j_2/10} + 2^{-(j-j_2)/10})
\|f_{j_1}\|_{Z} \|g_{j_2}\|_{Z}
\end{equation}
whenever $|j_1-j| \leq 10$ and $j_2 \leq j+11$.

We shall prove \eqref{high-low-targ} in Section \ref{high-low-sec}.  We leave this for now and turn to the high-high estimate \eqref{hhest}.
Here we cannot afford to discard the weights $w$.  By the triangle inequality in $l^2$, we can bound the left-hand side of \eqref{hhest} by
$$
\sum_{j_1,j_2:|j_1-j_2| \leq 1} (\sum_{j < j_1-10,j_2-10} \| 1_{A_j} \frac{w}{\langle \tau - \xi^2 \rangle} (\frac{f_{j_1}}{w} * \frac{g_{j_2}}{w}) \|_{Z}^2)^{1/2}.$$
By Schur's test again, it thus suffices to show that
$$
(\sum_{j < j_1-10,j_2-10} \| 1_{A_j} \frac{w}{\langle \tau - \xi^2 \rangle} (\frac{f_{j_1}}{w} * \frac{g_{j_2}}{w}) \|_{Z}^2)^{1/2}
\leq C \|f_{j_1}\|_{Z} \|g_{j_2}\|_{Z}$$
whenever $|j_1-j_2| \leq 1$.  Applying \eqref{localization} once more, and estimating the $Z$ norm by the $Y$ norm\footnote{This reflects the fact that it is very difficult for the high-high interaction to return to the parabola $\tau = \xi^2$, especially given our use of the weight $w$ to localize to the upper half-plane $\tau > 0$.},
we can simplify this a little as
\begin{equation}\label{high-high-targ}
 \| 1_{A_{\leq j_1-9}} \frac{w}{\langle \tau - \xi^2 \rangle} (\frac{f_{j_1}}{w} * \frac{g_{j_2}}{w}) \|_{Y}
\leq C \|f_{j_1}\|_{Z} \|g_{j_2}\|_{Z}.
\end{equation}
This estimate shall be proven in Section \ref{high-high-sec}.

Thus to conclude the proof of Proposition \ref{construct} (and hence Proposition \ref{func-sec}) it suffices to prove \eqref{high-low-targ}
and \eqref{high-high-targ}.  In many cases (basically, when at least two of $f, g, f*g$ are far from the parabola),
these inequalities can be established through Young's inequality, Proposition \ref{space-est}, and the resonance estimate \eqref{resonance}.
However, when two of $f,g,f*g$ are close to the parabola we need a further (standard) bilinear estimate, to which we now turn.

\section{Bilinear estimates near the parabola}

We give a standard bilinear estimate.

\begin{proposition}[Bilinear estimate]\label{bil-halt}
  Let $f, g$ be test functions supported on $A_{j_1}$ and $A_{j_2}$ respectively.  Suppose also that there is $D \geq 0$ such that
$|\xi_1 - \xi_2| \geq D$ whenever $(\tau_1,\xi_1)$ lies in the support of $f$ and $(\tau_2,\xi_2)$ lies in the support of $g$
(this hypothesis is vacuous if $D=0$).
Then
$$ \| f * g \|_{L^2_\xi L^2_\tau} \leq C 2^{j_1+j_2} \langle D \rangle^{-1/2} \| f \|_{\hat X^{-1,1/2,1}} \| g \|_{\hat X^{-1,1/2,1}}.$$
\end{proposition}

\begin{proof} 
Let $f_{d_1}$ be the restriction of $f$ to $B_{d_1}$, and similarly let $g_{d_2}$ be the restriction of
$g$ to $B_{d_2}$.  By \eqref{fx} we have
$$ \|f\|_{\hat X^{-1,1/2,1}} = 2^{-j_1} \sum_{d_1 \geq 0} 2^{d_1/2} \| f_{d_1} \|_{L^2_{\tau_1} L^2_{\xi_1}}; \quad
\|g\|_{\hat X^{-1,1/2,1}} = 2^{-j_2} \sum_{d_2 \geq 0} 2^{d_2/2} \| g_{d_2} \|_{L^2_{\tau_2} L^2_{\xi_2}}
$$
and so by the triangle inequality it suffices to show that
\begin{equation}\label{fudgud}
 \| f_{d_1} * g_{d_2} \|_{L^2(\R \times \R)} \leq C 2^{(d_1+d_2)/2} (2^{d_1/2} + 2^{d_2/2} + D)^{-1/2}
\| f_{d_1} \|_{L^2_{\tau_1} L^2_{\xi_1}} \| g_{d_2} \|_{L^2_{\tau_2} L^2_{\xi_2}}
\end{equation}
for each $d_1, d_2 \geq 0$.  Note that this is a little stronger than we need, as
$(2^{d_1/2} + 2^{d_2/2} + D)^{-1/2}$ is better than $\langle D \rangle^{-1/2}$, but we shall use this improvement in Corollary \ref{bil-dual} below.

Fix $d_1,d_2$; we may take $d_1 \geq d_2$ by symmetry. From Cauchy-Schwarz we have
\begin{align*}
\| f_{d_1} * g_{d_2} \|_{L^2(\R \times \R)}^2 &=
\int_\R \int_\R (\int_{B_{d_1} \cap ((\tau,\xi) - B_{d_2}); |\xi_1-\xi_2| \geq D} f_{d_1}(\tau_1,\xi_1) g_{d_2}(\tau_2,\xi_2)\ d\tau_1 d\xi_1)^2 d\tau d\xi\\
&\leq
\int_\R \int_\R \int_\R \int_\R f_{d_1}(\tau_1,\xi_1)^2 g_{d_2}(\tau_2,\xi_2)^2\ d\tau_1 d\xi_1 d\tau d\xi \\
&\quad\quad \sup_{\tau,\xi} |\{ (\tau_1,\xi_1) \in B_{d_1} \cap ((\tau,\xi) - B_{d_2}): |\xi_1-\xi_2| \geq D \}| \\
&= \| f_{d_1} \|_{L^2_{\tau_1} L^2_{\xi_1}}^2 \| g_{d_2} \|_{L^2_{\tau_2} L^2_{\xi_2}}^2\\
&\quad\quad \sup_{\tau,\xi} |\{ (\tau_1,\xi_1) \in B_{d_1} \cap ((\tau,\xi) - B_{d_2}): |\xi_1-\xi_2| \geq D \}|,
\end{align*}
where we are using the convention \eqref{tau-convention}.  Thus it suffices to show that
$$ |\{ (\tau_1,\xi_1) \in B_{d_1} \cap ((\tau,\xi) - B_{d_2}): |\xi_1-\xi_2| \geq D \}| \leq C 2^{d_1+d_2} / (2^{d_1/2} + D).$$
Observe that if $(\tau_1,\xi_1)$ lies in the above set, then $\tau_1 = \xi_1^2 + O( 2^{d_1} )$,
$\tau_2 = \xi_2^2 + O( 2^{d_2} )$, and thus $\tau = \xi_1^2 + \xi_2^2 + O( 2^{d_1} )$.  From the parallelogram identity
$$ \xi_1^2 + \xi_2^2 = \frac{1}{2}( \xi^2 + (\xi_1-\xi_2)^2 )$$
we thus have
$$ (\xi_1 - \xi_2)^2 = 2\tau - \xi^2 + O( 2^{d_1} ).$$
On the other hand, we have $|\xi_1 - \xi_2| \geq D$.  Elementary algebra then shows that $\xi_1 - \xi_2$ is confined to
a set of measure at most $O( 2^{d_1} / (2^{d_1/2} + D) )$.  Thus $\xi_1$ is also confined to a set of similar measure.  For fixed $\xi_1$
and $\xi_2$, $\tau_2$ (and hence $\tau_1$) is confined to an interval of length $O( 2^{d_2})$.  The claim then follows from Fubini's theorem.
\end{proof}

We can dualize this to obtain

\begin{corollary}[Dual bilinear estimate]\label{bil-dual}  Let $D \geq 0$, and suppose $\Omega_1 \subseteq A_{j_1}$, $\Omega \subseteq A_j$
be regions such that $|\xi_1 + \xi| \geq D$ whenever $(\tau_1,\xi_1) \in \Omega_1$ and $(\tau,\xi) \in \Omega$.
Then for any $f$ supported in $\Omega_1$, any test function $g$, and any $d \geq 0$, we have
$$ 2^{-d/2} \| f*g \|_{L^2_\xi L^2_\tau(\Omega \cap B_d)} \leq C 2^{j_1} (2^{d/2}+ D)^{-1/2} \| f \|_{\hat X^{-1,1/2,1}} \|g\|_{L^2_\xi L^2_\tau}.$$
\end{corollary}

\begin{proof} We can take $f,g$ to be non-negative.  By duality we can write
$$ 2^{-d/2} \| f*g \|_{L^2(\Omega \cap B_d)} = \int_{\R \times \R} f*g(\tau,\xi) h(\tau,\xi)\ d\tau d\xi$$
for some non-negative $h$ supported in $\Omega \cap B_d$ with $\|h\|_{L^2(\Omega \cap B_d)} = 2^{-d/2}$.  
We can then use the Fubini-Tonelli theorem, Cauchy-Schwarz, and Fubini-Tonelli again to write
\begin{align*}
\int_{\R \times \R} f*g(\tau,\xi) h(\tau,\xi)\ d\tau d\xi &=
\int_{\R \times \R} (\int_{\R \times \R} f(\tau_1,\xi_1) h(\tau_1+\tau_2,\xi_1+\xi_2)\ d\tau_1 d\xi_1) g(\tau_2,\xi_2)\ d\tau_2 d\xi_2 \\
&\leq 
(\int_{\R \times \R} (\int_{\R \times \R} f(\tau_1,\xi_1) h(\tau_1+\tau_2,\xi_1+\xi_2)\ d\tau_1 d\xi_1)^2\ d\tau_2 d\xi_2)^{1/2}\\
&\quad \times \|g\|_{L^2_\tau L^2_\xi(\R \times \R)} \\
&= \| f_-*h \|_{L^2(\R \times \R)} \|g\|_{L^2_\tau L^2_\xi(\R \times \R)}
\end{align*}
where $f_-$ is the reflection of $f$.
On the other hand, by decomposing $f = \sum_{d_1} f_{d_1}$, where each $f_{d_1}$ is supported on $B_{d_1}$, and using \eqref{fudgud}, \eqref{fx}
we have
\begin{align*}
 \| f_-*h \|_{L^2(\R \times \R)} &\leq \sum_{d_1}
2^{(d_1+d)/2} (2^{d_1/2} + 2^{d/2} + D)^{-1/2}
\| f_{d_1} \|_{L^2_{\tau_1} L^2_{\xi_1}} 2^{-d/2}\\
&\leq (2^{d/2} + D)^{-1/2} 2^{j_1} \|f\|_X
\end{align*}
and the claim follows.
\end{proof}

\section{High-low interactions}\label{high-low-sec}

We now prove the high-low interaction estimate \eqref{high-low-targ}.  Recall that we have
$|j_1-j| \leq 10$, $j_2 \leq j+11$.

Let us first dispose of the easy case $j_2=0$.  In this case we use \eqref{x-trivial-right}
followed by Young's inequality and Proposition \ref{space-est} to estimate
\begin{align*}
\| \frac{1_{A_j}}{\langle \tau - \xi^2 \rangle} (f_{j_1} * g_0) \|_{\hat X^{-1,1/2,1}}
&\leq C 2^{-j_1} \| f_{j_1} * g_0 \|_{L^2_\xi L^2_\tau} \\
&\leq C 2^{-j_1} \| f_{j_1} \|_{L^2_\xi L^2_\tau} \|g_0\|_{L^1_\xi L^1_\tau} \\
&\leq C \| f_{j_1} \|_Z \|g_0 \|_Z
\end{align*} 
which is acceptable.  Thus we may restrict attention to the case $j_2 > 0$.
Applying the resonance estimate \eqref{resonance}, we obtain
$$
\max( \langle \tau - \xi^2 \rangle, \langle \tau_1 - \xi_1^2 \rangle, \langle \tau_2 - \xi_2^2 \rangle)
\geq 2^{-20} 2^{j+j_2}.$$
Thus we may restrict one of $f_{j_1}$, $g_{j_2}$, or $f_{j_1}*g_{j_2}$ to the region $B_{\geq j+j_2 - 20}$.

Let us first consider the case when the high-frequency input $f_{j_1}$ is restricted to the region $B_{\geq j+j_2 - 20}$.  We can split this
case into two sub-cases, depending on whether we measure $g_{j_2}$ using $\hat X^{-1,1/2,1}$ or using $Y$.  If we use $Y$, then we use
H\"older's inequality in $\tau$, followed by Young's inequality and Proposition \ref{space-est}, \eqref{fy} to conclude
\begin{align*}
 \| \frac{1_{A_j}}{\langle \tau - \xi^2 \rangle} (f_{j_1} * g_{j_2}) \|_{\hat X^{-1,1/2,1}}
&\leq C 2^{-j} \| f_{j_1} * g_{j_2} \|_{L^2_\xi L^\infty_\tau}\\
&\leq C 2^{-j} \|f_{j_1} \|_{L^2_\xi L^2_\tau} \| g_{j_2} \|_{L^1_\xi L^2_\tau} \\
&\leq C 2^{-j} 2^j 2^{-(j+j_2)/2} \|f_{j_1}\|_{Z} 2^{j_2/2} \|g_{j_2}\|_Y \\
&\leq C 2^{-j_2/10} \|f_{j_1}\|_{Z} \|g_{j_2}\|_Y
\end{align*}
which is acceptable.  If we instead measure $g_{j_2}$ using $\hat X^{-1,1/2,1}$, we decompose into the regions $B_d$ and use
Corollary \ref{bil-dual} with $D=0$ (and with $f_{j_1}$ and $g_{j_2}$ swapped) followed by Proposition \ref{space-est}
to estimate\footnote{Note that if $j_2$ is substantially smaller than $j$ then we can take $D$ as large
as $2^j$, which leads to much better estimates.  However the $j=j_2$ case contains the delicate ``parallel interaction'' in which no gain occurs.}
\begin{align*}
 \| \frac{1_{A_j}}{\langle \tau - \xi^2 \rangle} (f_{j_1} * g_{j_2}) \|_{\hat X^{-1,1/2,1}}
&\leq C \sum_{d} 2^{-j} 2^{-d/2} \| f_{j_1} * g_{j_2} \|_{L^2_\xi L^2_\tau(A_j \cap B_d)} \\
&\leq C \sum_{d} 2^{-j} 2^{j_2} 2^{-d/4} \|g_{j_2} \|_{\hat X^{-1,1/2,1}} \|f_{j_1} \|_{L^2_\tau L^2_\xi} \\
&\leq C 2^{-j} 2^{j_2} \|g_{j_2} \|_{\hat X^{-1,1/2,1}} 2^{j_1} 2^{-(j+j_2)/2} \|f_{j_1} \|_{Z} \\
&\leq 2^{-(j-j_2)/10} \|f_{j_1}\|_{Z} \|g_{j_2}\|_{\hat X^{-1,1/2,1}}
\end{align*}
which is acceptable.  Thus we may now restrict $f_{j_1}$ to the region $B_{< j+j_2-20}$.  From Lemma \ref{paste} we may now measure $f_{j_1}$ in 
$\hat X^{-1,1/2,1}$ instead of $Z$.

We next consider the case when $g_{j_2}$ is restricted to the region $B_{\geq j+j_2-20}$.  
We may assume that $j_2 \leq j-10$, since when $j-10 < j_2 \leq j+11$ the situation here is essentially identical to the previous case (but with the roles of $f_{j_1}$ and $g_{j_2}$ reversed).
We subdivide the domain $A_{j_1} \cap B_{<j+j_2-20}$ of 
$f_{j_1}$ into disjoint slabs, where on each slab the frequency variable $\xi_1$ is localized to an interval $I$ of length $2^{j_2}/100$, and write $f_{j_1} = \sum_I f_{j_1,I}$ accordingly.  Because $g_{j_2}$
is localized to $A_{j_2}$, we see that the functions $f_{j_1,I} * g_{j_2}$ have finite overlap in the $\xi$ variable. Thus by square-summing in $I$ it would suffice to establish the estimate
$$  \| \frac{1_{A_j}}{\langle \tau - \xi^2 \rangle} (f_{j_1,I} * g_{j_2}) \|_{\hat X^{-1,1/2,1}}
\leq C (2^{-j_2/10} + 2^{-(j-j_2)/10}) \|f_{j_1,I}\|_{\hat X^{-1,1/2,1}} \|g_{j_2}\|_{Z}.$$
But from dyadic decomposition into $B_d$ regions and Corollary \ref{bil-dual} (which now applies with $D = 2^{j_2}/10$, say) followed
by Proposition \ref{space-est} we have
\begin{align*}
\| \frac{1_{A_j}}{\langle \tau - \xi^2 \rangle} (f_{j_1,I} * g_{j_2}) \|_{\hat X^{-1,1/2,1}}
&\leq C \sum_{d} 2^{-j} 2^{-d/2} \| f_{j_1,I} * g_{j_2} \|_{L^2_\xi L^2_\tau(A_j \cap B_d)} \\
&\leq C \sum_{d} 2^{-j} 2^{j_1} (2^{d/2} + 2^{j_2})^{-1/2} \|f_{j_1,I}\|_{\hat X^{-1,1/2,1}} \| g_{j_2}\|_{L^2_\xi L^2_\tau} \\
&\leq C j_2 2^{-j_2/2} \|f_{j_1,I}\|_{\hat X^{-1,1/2,1}} \|g_{j_2}\|_{Z} \\
\end{align*}
which is acceptable.  

The only remaining case is when we restrict the output $f_{j_1} * g_{j_2}$ to the region $B_{\geq j+j_2-20}$.  By Proposition \ref{space-est}, it
then suffices to show that
$$ 2^{-j} 2^{-(j+j_2)/2} \| f_{j_1} * g_{j_2} \|_{L^2_\xi L^2_\tau}
\leq C (2^{-j_2/10} + 2^{-(j-j_2)/10}) \|f_{j_1}\|_{\hat X^{-1,1/2,1}} \|g_{j_2}\|_{Z}.$$
If we measure $g_{j_2}$ in $\hat X^{-1,1/2,1}$, then from Proposition \ref{bil-halt} we have
\begin{align*}
  2^{-j} 2^{-(j+j_2)/2} \| f_{j_1} * g_{j_2} \|_{L^2_\xi L^2_\tau}
&\leq C 2^{-j} 2^{-(j+j_2)/2} 2^{j_1} 2^{j_2} \|f_{j_1}\|_{\hat X^{-1,1/2,1}} \|g_{j_2}\|_{\hat X^{-1,1/2,1}} \\
&\leq C 2^{-(j-j_2)/10} \|f_{j_1}\|_{\hat X^{-1,1/2,1}} \|g_{j_2}\|_{\hat X^{-1,1/2,1}}
\end{align*}
which is acceptable.  On the other hand, if we measure $g_{j_2}$ in $Y$, then from Young's inequality, Proposition \ref{space-est}, and
\eqref{fy} we have
\begin{align*}
 2^{-j} 2^{-(j+j_2)/2} \| f_{j_1} * g_{j_2} \|_{L^2_\xi L^2_\tau}
&\leq C 2^{-j} 2^{-(j+j_2)/2} \|f_{j_1} \|_{L^2_\xi L^1_\tau} \|g_{j_2}\|_{L^1_\xi L^2_\tau} \\
&\leq C 2^{-j} 2^{-(j+j_2)/2} 2^{j_1} \|f_{j_1} \|_{\hat X^{-1,1/2,1}} 2^{j_2/2} \|g_{j_2}\|_Y \\
&\leq C 2^{-j_2/10} \|f_{j_1} \|_{\hat X^{-1,1/2,1}}
\end{align*}
which is also acceptable.  This concludes the proof of \eqref{high-low-targ}.

\section{High-high interactions}\label{high-high-sec}

We now prove the high-high estimate \eqref{high-high-targ}.  Recall that $|j_1-j_2| \leq 1$
and we are operating under the assumption \eqref{assumption}.  We may also assume $j_1 \geq 9$ since the claim is
vacuous otherwise.  We can assume that $f_{j_1}$ is supported on one half-space, say $\{\xi > 0\}$, and that $g_{j_2}$ is supported on the other
half-space $\{\xi < 0\}$, since if they are both supported the same half-space then their convolution will not intersect
$A_{\leq j_1-9}$.  In particular we can ensure that the $\xi$-supports of $f_{j_1}$ and $g_{j_2}$ are separated by at least $2^{j_1}/10$.

We need some preliminary convolution estimates.
From Proposition \ref{bil-halt} (with $D = 2^{j_1}/10$) we have
$$ \| f_{j_1} * g_{j_2} \|_{L^2_\xi L^2_\tau} \leq
C 2^{j_1} 2^{j_2} 2^{-j_1/2} \| f_{j_1} \|_{\hat X^{-1,1/2,1}} \| g_{j_2} \|_{\hat X^{-1,1/2,1}}$$
while from Young's inequality and Proposition \ref{space-est}
\begin{align*}
\| f_{j_1} * g_{j_2} \|_{L^2_\xi L^2_\tau} 
&\leq \|f_{j_1} \|_{L^2_\xi L^2_\tau} \|g_{j_2} \|_{L^1_\xi L^1_\tau} \\
&\leq \|f_{j_1} \|_Y 2^{3j_2/2} \|g_{j_2} \|_{Z}
\end{align*}
and similarly
$$ \| f_{j_1} * g_{j_2} \|_{L^2_\xi L^2_\tau} \leq 2^{-2j_1} 2^{3j_1/2} \|f_{j_1} \|_{Z}
\|g_{j_2}\|_Y.$$
Putting all these estimates together, we obtain
$$
 \| f_{j_1} * g_{j_2} \|_{L^2_\xi L^2_\tau}
\leq C 2^{2j_1} \|f_{j_1} \|_{Z} \|g_{j_2} \|_{Z}.$$
In a similar spirit, from H\"older's inequality, Young's inequality and Proposition \ref{space-est} we have
\begin{align*}
\| \langle \xi \rangle^{-1} f_{j_1} * g_{j_2} \|_{L^2_\xi L^1_\tau} 
&\leq C\| f_{j_1} * g_{j_2} \|_{L^\infty_\xi L^1_\tau} \\
&\leq C\| f_{j_1} \|_{L^2_{\xi} L^1_\tau} \| g_{j_2} \|_{L^2_\xi L^1_\tau} \\
&\leq C 2^{j_1} \| f_{j_1} \|_{Z} 2^{j_2} \| g_{j_2} \|_{Z};
\end{align*}
combining these estimates using \eqref{fy} we obtain 
\begin{equation}\label{fgj}
 \| f_{j_1} * g_{j_2} \|_Y \leq C 2^{2j_1}
\|f_{j_1} \|_{Z} \|g_{j_2} \|_{Z}
\end{equation}

We now return to \eqref{high-high-targ}.
First let us restrict $A_{\leq j_1-9}$ to the region $A_{\leq j_1-9} \cap B_{\geq 2j_1 - 10}$.  
In this case we discard the weights $w$ to obtain
\begin{equation}\label{abab}
 \| 1_{A_{\leq j_1-9} \cap B_{\geq 2j_1-10}} \frac{w}{\langle \tau - \xi^2 \rangle} (\frac{f_{j_1}}{w} * \frac{g_{j_2}}{w}) \|_{Y}
\leq C 2^{-2j_1} \| f_{j_1} * g_{j_2} \|_{Y},
\end{equation}
and \eqref{high-high-targ} in this case follows from \eqref{fgj}.

Thus it remains to consider the contribution in the domain $A_{\leq j_1-9} \cap B_{< 2j_1-10}$ (i.e. the contribution near
the frequency origin).  We now finally invoke \eqref{assumption}; we shall assume that $f_{j_1}$ lies in $B_{\leq 2j_1 - 100}$ since
the other case is almost identical (recall that $|j_1-j_2| \leq 1$).  In particular, we have $\tau_1 \geq 2^{2j_1}/10$ and $|\tau| \leq 2^{2j_1}/100$,
which forces $\tau_2 \leq - 2^{2j_1}/20$ by \eqref{tau-convention}.  Thus we now have a large weight on $g_{j_2}$: $w(\tau_2,\xi_2) \geq c 2^{20j_1}$.
On the other hand, we make the elementary observation that
$$ \frac{w(\tau,\xi)}{\langle \tau - \xi^2 \rangle} \leq C 2^{-2j_1} 2^{20j_1}$$
for $(\tau,\xi)$ in $A_{\leq j_1-9} \cap B_{< 2j_1-10}$.  Thus in this case we again have
\eqref{abab}, and again \eqref{high-high-targ} in this case follows from \eqref{fgj}.  This concludes the proof of \eqref{high-high-targ},
and thus of Proposition \ref{construct} and Proposition \ref{func-sec}.

\end{document}